\documentclass[manuscript, 11pt]{acmart}
\usepackage{lineno}

\renewcommand\footnotetextcopyrightpermission[1]{} 
\usepackage[utf8]{inputenc}

\title{Elitism in Mathematics and Inequality}
\author{Ho-Chun Herbert Chang}
\affiliation{%
    \institution{University of Southern California}
    \department{Annenberg School of Communication and Journalism}
    \department{Information Sciences Institute}
    }
    
\author{Feng Fu}
    \affiliation{%
    \institution{Dartmouth College}
    \department{Department of Mathematics}
    }
\date{February 2020}

\begin{document}
\begin{abstract}
\textbf{Abstract: }The Fields Medal, often referred as the Nobel Prize of mathematics, is awarded to no more than four mathematician under the age of 40, every four years. In recent years, its conferral has come under scrutiny of math historians, for rewarding the existing elite rather than its original goal of elevating mathematicians from under-represented communities~\cite{barany2018fields,barany2015myth}. Prior studies of elitism focus on citational practices~\cite{hirsch2005index} and sub-fields~\cite{rossi2017genealogical,gargiulo2016classical}; the structural forces that prevent equitable access remain unclear. 
Here we show the flow of elite mathematicians between countries and lingo-ethnic identity, using network analysis and natural language processing on 240,000 mathematicians and their advisor-advisee relationships. We found that the Fields Medal helped integrate Japan after WWII, through analysis of the elite circle formed around Fields Medalists.  Arabic, African, and East Asian identities remain under-represented at the elite level. Through analysis of inflow and outflow, we rebuts the myth that minority communities create their own barriers to entry. Our results demonstrate concerted efforts by international academic committees, such as prize giving, are a powerful force to give equal access. We anticipate our methodology of academic genealogical analysis can serve as a useful diagnostic for equality within academic fields.
\end{abstract}

\maketitle

Although mathematics is often framed as objective and egalitarian, its access is not equally conferred. Recent attention has been given to the Fields Medal, one of the most prestigious awards in math, and its elite community. 
When the award was first conceived in 1930, it was in part designed to assuage international tensions~\cite{barany2018fields}. The award was intentionally given to individuals that would otherwise not receive any recognition, rather than the best young mathematician.

Using social network analysis (SNA) and neural-based natural language processing (NLP), this paper analyses the flow of elite mathematicians between nations and lingo-ethnic categories.  Analysis was performed on the Mathematics Genealogy Project, one of the most complete advisor-advisee databases maintained today with more than 240,000 mathematicians.

Results demonstrates the self-reinforcing behavior among the elite level in mathematics. This contrasts with prior conferral of the Fields Medal, which was a positive force in mending international relations, such as integrating Japan and Germany after World War II~\cite{parshall2009internationalization}. We propose the Fields Medal can be used today to improve accessibility of mathematics to minority groups. 

The classifier for lingo-ethnic identity is textual, it would be more accurate to say we classify specific languages that overlap significantly with ethnic or cultural identity. While the use of lingo-ethnic categorization as identity is shallow, our principal aim is to show, even at the most basic definitions of ethnicity or culture through language, we find evidence of inequality. This paper also offers a methodological contribution. We show that combining network analysis, neural-based natural language processing (NLP), and well-maintained academic databases can serve as a powerful diagnosis for access and equity, and improve the practice of science.

Several studies on elitism within the production of mathematical knowledge have been conducted. Methods draw predominantly from the complex network perspective~\cite{zeng2017science}, leveraging network repositories such as citation and bibliometric networks. Gargiulo et al. studied the entire, connected giant component of the mathematical genealogy project, enriching the data using data mining techniques~\cite{gargiulo2016classical}. They work focused on integrating math history with temporal network analysis, noting the fields evolution based on country, discipline, and the structure of scientific families. 

Prior investigated about the relationship between scientific mentorship and winning the Fields Medal or Wolf Prize, but results were inconclusive. Rossi et al. studied the role of advisor-advisee relationships ~\cite{rossi2017genealogical}. 
They propose the \textit{genealogy index}, adapted from the \textit{h-index} which was initially developed by Hirsh~\cite{hirsch2005index}. 
 Malmgren et al. studied the role of mentorship on protégé performance, focused on metrics of academic success like publication record~\cite{malmgren2010role}. Beyond scholarship, studies have also considered hiring practices~\cite{clauset2015systematic} and departmental prestige~\cite{myers2011mathematical}.

The lack of metadata in these genealogies has limited the scope of investigation. 
This paper places elite community network flow as the focal point, contrasting the historical focus on the nation-state with the modern focus of identity.

\section*{Historical Networks of Elite Migration}
\begin{figure}[!htb]
    \centering
    \includegraphics[width = 1.0\linewidth]{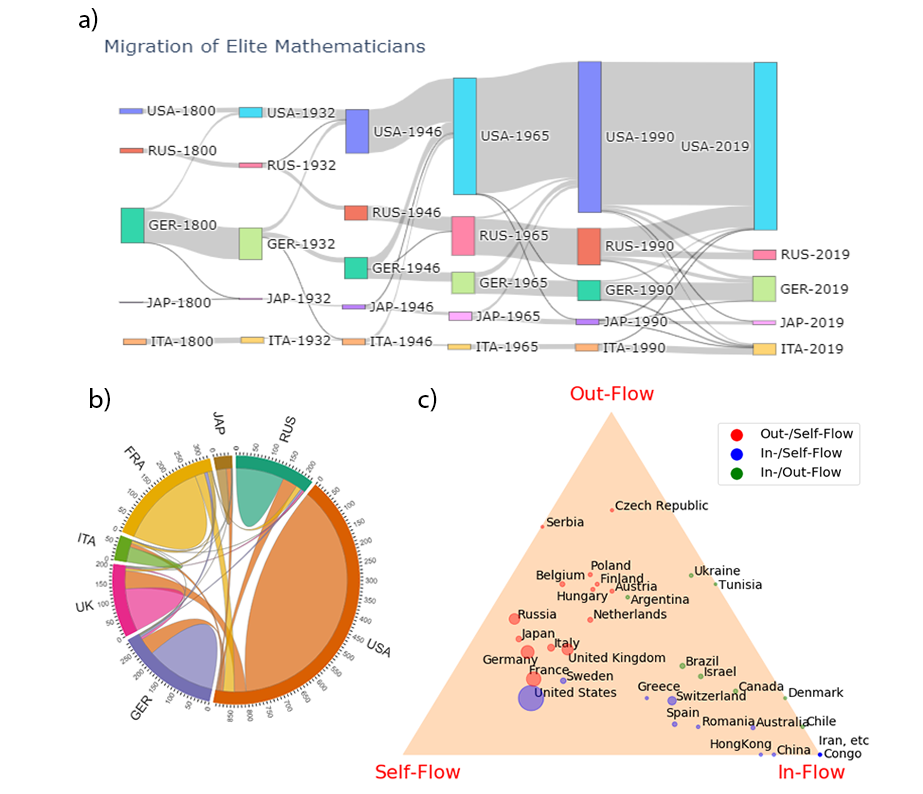}
    \caption{\textbf{Mobility patterns of mathematicians among countries (traditionally strong en mass)}. \ref{fig:countries}a), the migration of elite mathematicians from 1800 to the present. \ref{fig:countries}b), the net flow of elite mathematicians between 7 key countries. \ref{fig:countries}c), the flow analysis of elite mathematicians between countries. Exporting means mathematicians flow out a country, importing means they flow in, and selfish means they are retained. Many countries only import at the elite level (\textit{etc} in the bottom right). A full list is available in the appendix.}
    \label{fig:countries}
\end{figure}

We begin with a sketch of history. Figure 1a) captures the migration of elite mathematicians between five key countries. The subgroup of elites was created by aggregating the shortest paths between Fields Medalists. This ensures that the full graph is connected, and conceptually, denotes a minimal graph that connects all the medalists together. Here, migration is determined by comparing where a mathematician earned their Ph.D. and where their students earned their Ph.D. It is reasonable to assume primary advisors have moved to the same country as their advisees.

Prior to WWII, Western European countries were the strong-holds of mathematical thought. Notably, France and Germany contained the highest proportion of elite mathematicians. Many Japanese mathematicians studied in Germany, before returning to Japan, as part of modernization during the Meiji restoration. Examples inclue Rikitaro Fujisawa, who studied at the Unviersity of Strasbourg with Elwin Christofeel, before returning~\cite{chikara2013intersection}. He was instrumental to reforming mathematics education in Japan.

The flow chart reveals mass flows of researchers due to historical events. By 1932, the Holocaust led to mass migration from Germany to the United States and other European countries, which accounts for the drop in green volume, including prominent scientist Albert Einstein. Similarly, we observe large amounts of outflow from Russia after the cold war, greatly diminishing the presence of Russia mathematicians after the 1990s, and the second Italian mass diaspora after WWII.
Beyond forced immigration, flow analysis also reveals the movement of reintegration. Japanese mathematicians immigrated to the United States following WWII, and continued throughout the 60s to the 90s. Twenty years later, Japanese mathematicians flowed back toward Japan. 

France is not shown in the Sankey flow chart (1a), but is historically one of the countries that produces the most elite mathematicians. The chord graph in 1b) shows the net flow of mathematicians over all time, with the color of the chord indicating net exports. The USA-GER chord is orange, which indicates a net outflow from USA to Germany. Only France exports more American mathematicians than it imports from the USA. In all other cases, the USA exports more to other countries.

Figure 1c) shows the flow dynamics on a country level. In-flow is defined as the number of incoming edges, out-flow as the number of outgoing edge, and self-flow the number of loops. These results are similar to Gargiulo et al.~ \cite{gargiulo2016classical} with two striking differences. First, the United States is a selfish and importing country at the elite level, whereas in general it is selfish and exporting. Secondly, there are many more importing countries compared to the general case, where most countries are exporting and selfish. Notice, many of the countries that are exporting and selfish are Western or part of the Soviet Union, where there were strong programs in mathematics. Other countries appear to import more at the elite level, because their "exports" are not as competitive as mathematicians exported from other countries.

These two points allow us to infer three things. First, elite mathematicians have more mobility, and in many cases can begin work in foreign countries. Second, the United States imports more compared to the general case, attracting more elite members. Third, countries considered traditional mathematics strong-holds can be observed in the lower left corner.

What this analysis tells us, beyond an exposé of diasporic history, is the fields medal served as a way to mediate tensions. In a similar way that the Olympics was held in Rome, Berlin, and Tokyo, the inclusion of internationally marginalized nations.

\section*{The Flow of Marginalized Identities   }
Upon analyzing the history of elite communities in mathematics, we turn to the present. As 1a) shows, today, there is significant flow between countries. lingo-ethnic categories of identity serve as a useful construct for understanding network flow. Figure 2a) shows the representation of identities, within three subgroups: all mathematicians (blue), mathematicians within the medalist subgroup, (green) and the medalists themselves (red).

\begin{figure}[!htb]
    \centering
    \includegraphics[width = 0.9\linewidth]{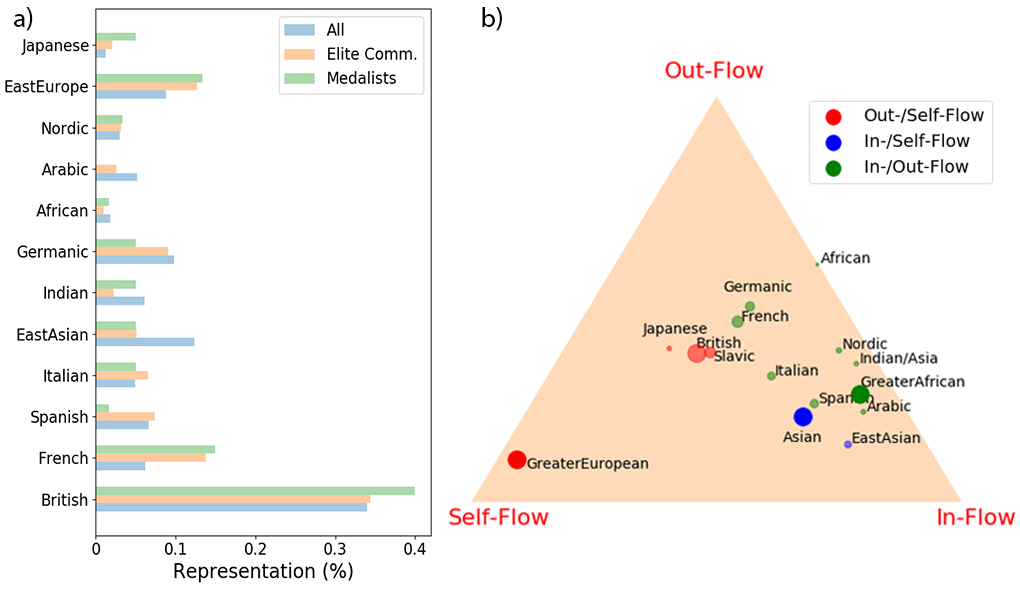}
    \caption{\textbf{Elite recognition and representation among lingo-ethnic identities}. \ref{fig:ethnicity}a), the proportion of mathematicians within the elite community by identity. The three subgroups are:  three subgroups: all mathematicians (blue), mathematicians within the medalist subgroup, (green) and the medalists themselves (red). \ref{fig:ethnicity}b), flow analysis by identity.}
    \label{fig:ethnicity}
\end{figure}

Fig.~\ref{fig:ethnicity} compares elite representation of subgroups relative to their actual proportions. For instance, there is a higher proportion of French medalists (14\%) compared to the general proportion (8\%). In contrast, there is a significant number of East Asian mathematicians (14\%) but very low representation in both the medalist family and medalists themselves (5\% each).

Upward sloping bars (left to right) mean medalists and medalist families are \textit{over-represented}; downward sloping bars indicate \textit{under-representation}.  Over-represented groups include British, French, Japanese, East European, and Nordic names. Underrepresented groups include East Asian and Germanic. Mathematicians with Arabic names are non-existent in Medalists and underrepresented in the elite community.

On the level of flow, Figure 2b) characterizes identities in terms of in-flow, out-flow, and self-flow. High in-flow means a higher likelihood of being mentored. High out-flow then corresponds to a greater likelihood to mentor others. High self-flow means higher likelihood of mentoring your own identity. The identity with the most self-flow is Japanese. However, when all mathematicians are considered, the Japanese are shown as green, that is to say opposite of selfish. This indicates reinforcing behavior only occurs at elite levels. 

However, once these groups are aggregated into larger groups--- Greater European, Asian, African and Arabic--- then differences become evident. European names has high self-reinforcing behavior, whereas Asians names and African and Arabic names are much lower in the number of self-loops. This dispels a common myth that minority groups, due to homophily, tend to group together. This myth insinuate that barriers to entry are self-inflicted. However, as we see from 2b), most minority groups are far away from the selfish pole, with a healthy balance of in-flow and out-flow. Rather, increases in the quantity of self-loops occurs in the greater European subgroup.

\section*{Old Strongholds, New Possibilities}
It is understandable that, when considering all mathematicians, that there is a high levels of self-flow--- studying in elite and often foreign institutions is a privilege. However, the fact that high self-flow in identity at the \textit{elite level} suggests institutions can do more to open access, given their greater access to resources. This has been the case for Japan.

Japan is unique among Asian countries and identities in that there are many Japanese Fields Medalists (3), with high representation in elite levels. Japan has been known for its rapid westernization during the Meiji restoration relative to other Asian counterparts. As early as 1872, their traditional form of math \textit{wasan} was replaced by western science. Prussia, rather than the United Kingdom, was the primary source of westernization, and led directly to the establishment of the University of Tokyo~\cite{parshall2009internationalization}. 

\begin{figure}[!htb]
    \centering
    \includegraphics[width = 1.0\linewidth]{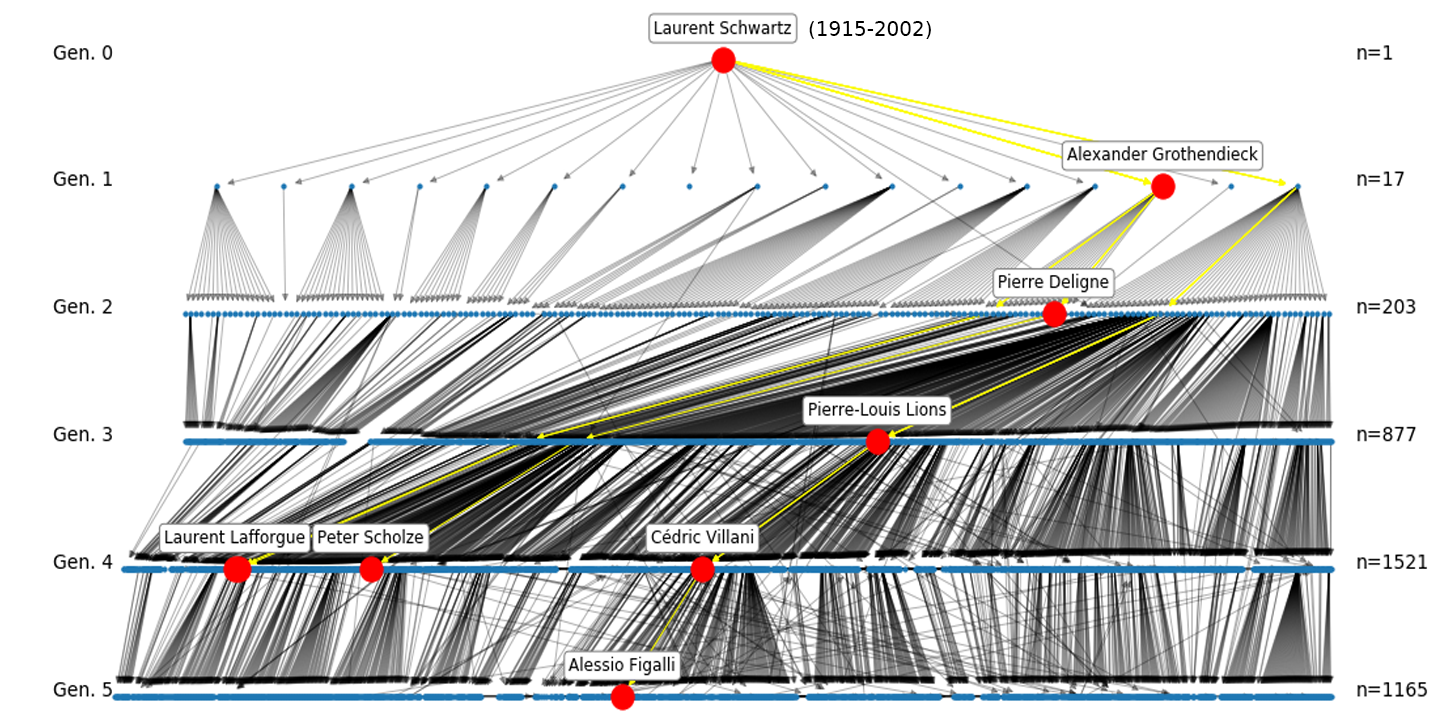}
    \caption{Super "reproductive" lineage of Fields Medalist Laurent Schwartz, up to five generations.}
    \label{fig:schwartz}
\end{figure}

After WWII, mathematicians sought to re-establish international ties and formed the International Congress of Mathematicians and a new International Mathematics Union (IMU). Marshall Stone, a proponent of this movement, said it clearly: "in considering American adherence to a Union, it must be borne in mind that we want nothing to do with an arrangement which excludes Germans and Japanese as such." Indeed, we find the ten founding members well-represented in the ternary diagrams, and not long after founding, the Soviet Union joined. Revisiting Fig. 1a), we discover the density of elite mathematicians in Japan increases after 1945. 

What this says, is the Fields Medal can improve the status of marginalized populations. Mathematics historian Barany captures this aspiration, believing the fields medal should help "sculpt the future, rather than reward the past~\cite{barany2018fields}." What we observe is the opposite, where the elite perpetuate the elite. Fig. 3 demonstrates this clearly, showing French Fields Medalist Laurent Schwartz and his lineage.

\begin{figure}[!htb]
    \centering
    \includegraphics[width = 1.0\linewidth]{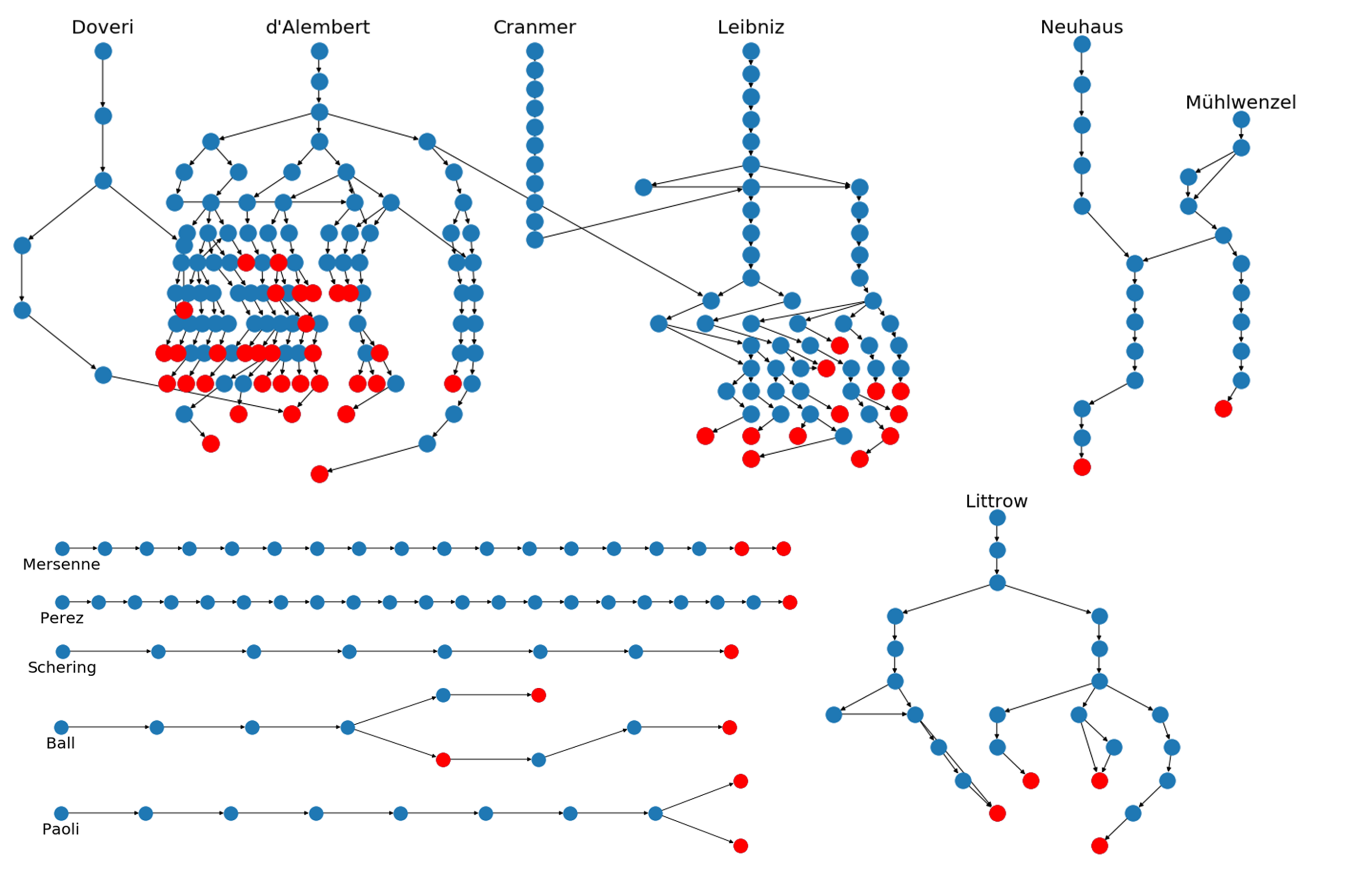}
    \caption{Ancestral Tree by Principle Advisor for each Fields Medalist. Each red dot denotes a Fields Medalist.}
    \label{fig:all-trees}
\end{figure}

Within 5 generations after Schwartz, 7 Fields Medalists emerge. In particular, Schwartz-Grothendieck-Deligne form a direct chain, as do Lions-Villani-Figalli. Note, Lions’ father Jacque -Louis Lions was also a student of Schwartz. In other words, 13.3\% of all Fields Medalists descended directly from Schwartz. Broadly, each of these all made contributions to some form of algebraic geometry or functional analysis. Fig.~\ref{fig:all-trees} further shows that all medalists can be traced to 9 connected components, with the largest one holding 44 out of 60 listed Fields Medalists.

\begin{figure}[!htb]
    \centering
    \includegraphics[width = 0.8\linewidth]{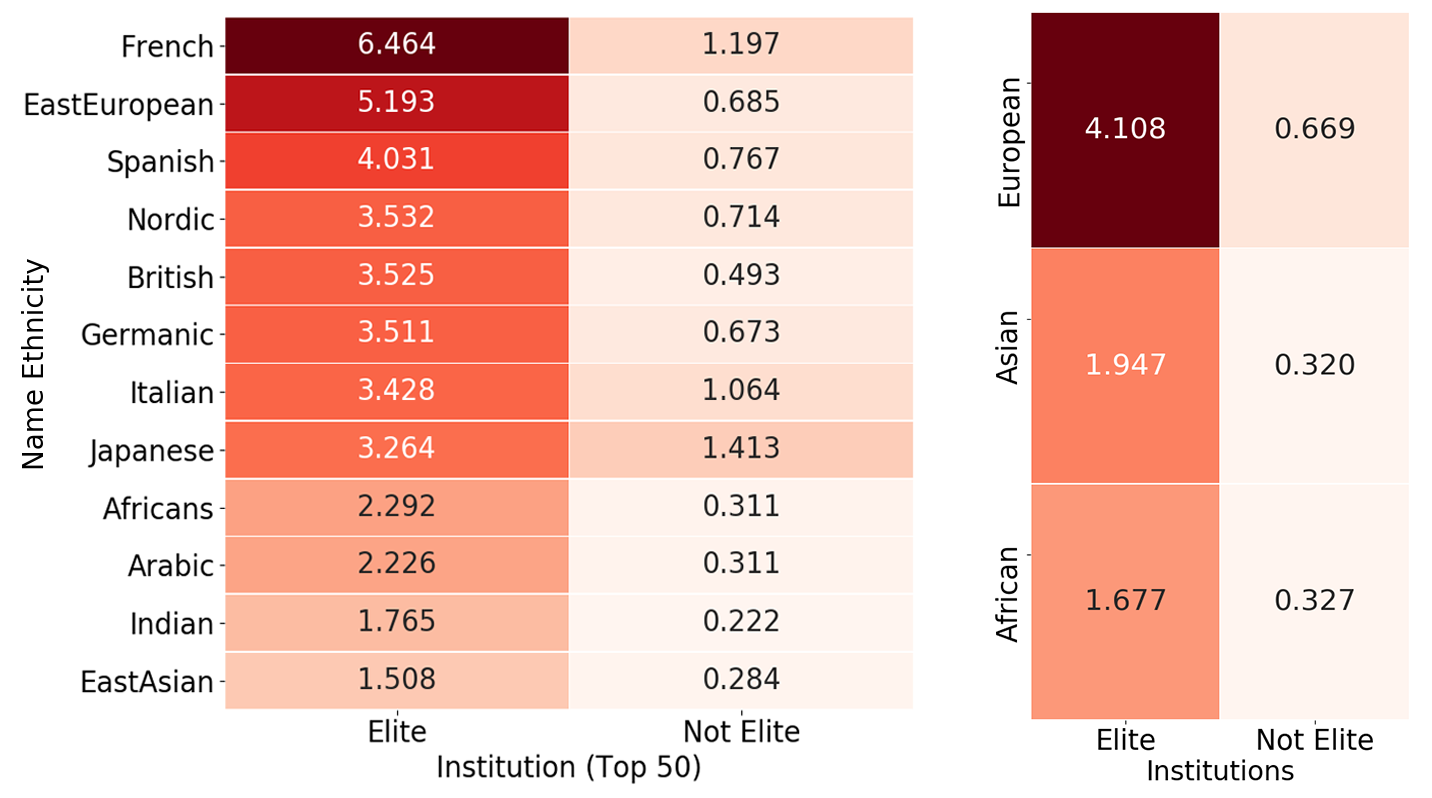}
    \caption{\textbf{The power ratio by identity.} Elite circles is a strong indicator of recognition. The power ratio is the likelihood of being part of the elite community, divided by the average likelihood. It is a measure of relative likelihood.}
    \label{fig:heat}
\end{figure}

These observations are not meant to diminish the achievements of great mathematicians. They do however show the Fields Medal has deviated from its commitment to elevate under-represented mathematicians.  Fig.~\ref{fig:heat} shows this succinctly in a tabular heatmap, which shows the power ratios. The power ratio (defined in Equation~\ref{eq:power-ratio}) is the conditional likelihood of being in the Fields Medalist Subgroup over the average probability of being in the group ($P = 0.00759$).

\begin{equation}\label{eq:power-ratio}
   PR =  \frac{P(\textit{Fields} ~|~ \textit{Institute \& Identity} )}{P(\textit{Fields})}
\end{equation}

A mathematician that is French and attends a Top 50 institution means they are 6.4 times more likely to gain membership into the elite circle. Here, the top 50 is defined as the top institutions attended by those in the elite group.
Note, we defined our Fields Medalist subgroup minimally, such that any other definition of subgroup would yield a higher power ratio. On the other hand, being East Asian and attending a Top 50 institution only affords you 1.4 times the likelihood of gaining membership into this elite circle.

From this diagram, we infer that institution plays a large role in elite membership. However, notice an East Asian mathematician a top 50 school is 4.5 times less likely to be included than a French mathematician attending a top 50 school. An Indian mathematician educated outside top 50 schools are 6 times less likely to be included than a French mathematician with the same education. Amongst non-elite institutions, being Japanese gives the best chance of inclusion, an after-effect of the efforts by the IMU.

\section*{Conclusion}

In 2014, the late Iranian mathematician Maryam Mirzakhani won the Fields Medal. A talented star herself, her groundbreaking work on dynamics and geometry was encouraged by her Ph.D. advisor Curtis McMullen, also a Fields Medalist, at the elite institution Harvard University. This is by no means downplaying her achievements; rather, it serves to show the power recognition and elite communities have---all of which membership she rightly earned. Although the Fields Medal should serve to recognize under-represented researchers, the proper cultivation of talent through mentorship and institutional support should be the starting point.

In our evaluation of the present, there is a large under-representation of minority groups in not just Field Medalists, but also in the elite circle for mathematics. While institutional prestige a big factor, lingo-ethnic identity is also found to be highly relevant, the widest gap being 4.5 times the power ratio even at elite institutions. Given that elite institutions have more resources, they can take a bigger role in generating higher access for marginalized groups. Flow analysis also dispels the myth that under-representation arises from homophily-driven self-selection.

Although the French stronghold shows the old forces that govern mathematical knowledge remain strong, the presence of Japanese scholars also shows concerted effort can be used as an integrating force. Concerted efforts by international academic committees, such as prize giving, are a powerful force to confer equal rights for knowledge production to traditionally marginalized groups. Beyond analysis, this network analytical methodology is a call for scientific communities to use advisor-advisee databases to open knowledge production and scientific access.

\bibliographystyle{unsrt}
\bibliography{references}

\section*{Methods}
\subsection*{Graph Construction}
The graph was constructed using the Mathematics Genealogy database. Nodes are mathematicians, and directed edges represent advisor-advisee relationships. The data set contained information (listed in order of completeness) on the academic, advisor-advisee links, school, PhD graduation year, country, and dissertation title and topic. The ID's of medalists were identified, then the shortest path was computed in a pairwise fashion. Analysis was conducted primarily using the Networkx package~\cite{hagberg2008exploring}.

The subgroup of elites was created by taking the union of shortest paths between Fields Medalists. Then, the full graph is connected, and denotes some form of minimal graph that connects all the medalists together. While it is possible to produce a minimal spanning tree, given the forest like structure of the genealogy, the shortest paths has more interpretive value.

\subsection*{Identity Classifier}
Since lingo-ethnic identity is not included in the Mathematics Genealogy Project, a separate classifier is required. The identity categories were labeled using the \textit{ethnicolr} package, which is a long-short term neural network (LSTM) trained on Wikipedia and the census~\cite{sood2018predicting}. Specifically, the LSTM was based off the seminal work of Graves and Schmidhuber~\cite{graves2005framewise}. This package has found use in evaluating under-representation in other STEM fields such as biomedicine~\cite{marschke2018last}. It achieves between 78\% to 81\% accuracy.  Potential shortcomings of neural methods for categorization is the accuracy levels. However, for 13 individual categories (which would result in 7.7\% accuracy if truly random), 81\% is quite high. Additionally, since we are interested in comparison within individual demographics, any bias would be carried forward since the group of all mathematicians supersets the medalist subgroup and medalists themselves. 
The goal of using this classifier is not to flatten definitions of identity, but to use the best available tools for inference, in absence of concrete data.

\subsection*{Flow Analysis}
Meso-graphs were constructed on attributes of each mathematician. To turn attributes into nodes, we constructed a mapping from mathematician to the meso-categories (lingo-ethnic and nationality of doctoral degree). Edges between meso-categories were simply the original directed-edges between mathematicians. Each edge is then weighed by the number of advisor-advisee relations between meso-categories.

\subsubsection*{Constructing Ternary Diagrams}
We constructed the ternary diagrams through analysis of the meso-network. Every meso-network can be represented by a its adjacency matrix, which we denote $M$. The diagonal then accounts for self-loops, the rows excluding the diagonal elements the out going edges, and the columns excluding the diagonal element the incoming edges. 
Explicitly, for meso-category indexed by $i$, we have the following definitions for in-flow (IF), out-flow (OF), and self-flow (SF).
$$
\begin{aligned}
SF_i &= M_{i,i} \\
IF_i &= \sum_{j \neq i} M_{i,j} \\
OF_i &= \sum_{j \neq i} M_{j,i}
\end{aligned}
$$
We then normalize these values to represent each meso-category as a point in three dimensional space.
$$
P_i = \Big( \frac{IF_i}{K_i}, \frac{OF_i}{K_i} , \frac{SF_i}{K_i}   \Big) \in [0,1]^3 \qquad 
\text{with } K_i = IF_i + OF_i + SF_i
$$
Note, all points lie on the plane described by $x + y + z = 1$. We then transform this planar section onto the two dimensional plane using a translation and two rotations.

\begin{equation}
\begin{aligned}
    P_i'  &= P_i - (0,0,1) \\
    P_i'' &= R_2 \circ R_1 (P_i')
\end{aligned}
\end{equation}
where $R_1$ rotates the plane up to the XY-plane, and $R_2$ aligns the simplex to the x-axis.

\section*{Acknowledgements}
  Both authors thanks the Mathematics Genealogy Project for generously providing data from its database for use in this research. H.C gratefully acknowledges the Annenberg Fellowship from the Annenberg School of Communication and Journalism, USC. F.F. gratefully acknowledges the Dartmouth Faculty Startup Fund, the Neukom CompX Faculty Grant, Walter \& Constance Burke Research Initiation Award and NIH Roybal Center Pilot Grant.


\end{document}